\documentclass[11pt]{article}
\usepackage[margin=1in]{geometry}
\usepackage{amsthm,amsmath}
\usepackage{latexsym,amssymb}
\usepackage{graphicx,float,color,fancybox,shapepar,setspace,hyperref}
\usepackage{subfigure}
\usepackage{pgf,tikz}
\usepackage[affil-it]{authblk}
\usepackage{indentfirst}
\usepackage[section]{placeins}
\hypersetup
{
colorlinks=true,
linkcolor=blue,
filecolor=blue,
urlcolor=blue,
citecolor=cyan,
}
\usetikzlibrary{arrows}
\newtheorem{thm}{Theorem}
\newtheorem{lem}{Lemma}
\newtheorem{cor}{Corollary}
\newtheorem{rk}{Remark}

\newcounter{remarkcounter}

 \def\qed{\hfill\square}

\def\~{\sim}

\def\qed{ \hfill $\blacksquare$}

\begin{document}

\title{Minimum degree and sparse connected spanning subgraphs}
\author{Ting HUANG$^1$, Yanbo ZHANG$^2$, Yaojun CHEN$^{1,}$\footnote{Corresponding author. Email: yaojunc@nju.edu.cn}\\
{\small $^1$School of Mathematics, Nanjing University, Nanjing 210093, China}\\
{\small $^2$School of Mathematical Sciences, Hebei Normal University, Shijiazhuang 050024, China}}
 \date{}
\maketitle

\begin{abstract}

Let $G$ be a connected graph on $n$ vertices and at most $n(1+\epsilon)$ edges with bounded maximum degree, and $F$ a graph on $n$ vertices with minimum degree at least $n-k$, where $\epsilon$ is a constant depending on $k$.  In this paper, we prove that $F$ contains $G$ as a spanning subgraph provided $n\ge 6k^3$, by establishing tight bounds for the Ramsey number $r(G,K_{1,k})$, where $K_{1,k}$ is a star on $k+1$ vertices. Our result generalizes and refines the work of Erd\H{o}s, Faudree, Rousseau, and Schelp (JCT-B, 1982), who established the corresponding result for $G$ being a tree. Moreover, the tight bound for $r(G,tK_{1,k})$ is also obtained.

\vskip 2mm
\noindent{\bf Keywords}: Spanning subgraph, Ramsey number, sparse graph, star
\end{abstract}

\section{Introduction}
  \noindent 
Let $G$ be a simple graph with vertex set $V(G)$ and edge set $E(G)$. For any $v\in V(G)$, the neighborhood of $v$ is denoted by $N(v)$ and $N[v]=N(v)\cup \{v\}$. 
  For $S\subseteq V(G)$, $G-S$ is the subgraph obtained from $G$ by deleting all vertices of $S$. Set $G_v=G-N[v]$.  The degree of $v$ is $d(v)=|N(v)|$. The minimum and maximum degrees of $G$ are denoted by $\delta(G)$ and $\Delta(G)$, respectively. The independence number of $G$ is denoted by $\alpha(G)$.
Define 
$$\alpha'(G) = \min\{\alpha(G_v):v\in V(G)\}.$$
Given a positive integer $n$ and a graph $F$, how can we ensure that $F$ contains all trees $T$ with $n$ vertices?  Except for the well-known folklore result that  $\delta(F)\geq n-1$ can guarantee this, Faudree, Rousseau, Schelp, and Schuster started to study this issue under the assumption that $\delta(F)$ is less than $n-1$ and  $\Delta(F)=n-1$ (that is necessary for $F$ containing all spanning trees), 
and proved the following result as early as 1980.
\begin{thm} (Faudree, Rousseau, Schelp, and Schuster \cite{FRSS})
    If $k \geq 3$ and $n \geq 3k^2-9k+8$, then every graph $F$ of 
order $n$ satisfying $\Delta(F) = n- 1$ and $\delta(F) \geq n -k$ contains all trees with $n$ vertices.
\end{thm}

Moreover, Erd\H{o}s, Faudree, Rousseau, and Schelp showed that if we do not require $\Delta(F)=n-1$, then the condition $\delta(F)\ge n-k$ can guarantee that $F$ contains a certain family of spanning trees, and established the following.   

\begin{thm} (Erd\H{o}s, Faudree, Rousseau, and Schelp \cite{EFRS}) \label{tree}
    If $k\geq 2$ and $n \geq 2(3k-2)(2k-3)(k- 2) + 1$, then every graph $F$ of order $n$ with $\delta(F) \geq n -k$ contains every tree $T$ of order $n$ and $\Delta(T)\leq n-2k+2$. 
\end{thm}
Note that the condition $\delta(F) \geq n -k$ is equivalent to the complement graph $\overline{F}$ of $F$ contains no star $K_{1,k}$,
Erd\H{o}s et al. also proposed that this sufficient condition is used to calculate some of the \emph{Ramsey numbers} for the pair tree-star. 
For two given graphs $G$ and $H$, the \emph{Ramsey number} $r(G,H)$ is defined as the smallest positive integer $N$ such that, for each graph $F$ on $N$ vertices, either $G$ is a subgraph
 of $F$ or $H$ is a subgraph of $\overline{F}$. Equivalently, it is also defined as the smallest positive integer $N$ such
that, for every red-blue edge-coloring of the complete graph $K_N$, there exists either a red subgraph isomorphic to $G$ or a blue subgraph isomorphic to $H$.

 
 Using Theorem \ref{tree}, Erd\H{o}s, Faudree, Rousseau, and Schelp obtained the following.


\begin{thm} (Erd\H{o}s, Faudree, Rousseau, and Schelp \cite{EFRS}) \label{tree star}
    Let $k\geq 2$ and $n \geq 2(3k-2)(2k-3)(k- 2) + 1$. If $T$ is a tree with $n$ vertices and $\alpha'(T)=\alpha'$, then 
    $$\max \{n, n + k - 1- \alpha'-\beta \} \leq r(T, K_{1,k}) \leq \max\{n, n + k - 1 - \alpha'\},$$ 
    where $\beta=0$ if $k\,|\, n+k-2-\alpha'$, and $\beta=1$ otherwise.  
\end{thm}

\begin{rk} \label{remark1} It is worth noting that the parameter $\alpha'(T)$ is closely related to $\Delta(T)$. If $\Delta(T)\leq n-2k+2$, then $\alpha'(T)\geq k-1$. This is because $|T_v|\geq 2k-3$ for all $v\in V(T)$. Since $T_v$ is also a bipartite graph, there must exist one independent set containing at least $k-1$ vertices in $T_v$, which implies that $\alpha(T_v)\geq k-1$ for all $v\in V(T)$. Hence we have $\alpha'(T)\geq k-1$. 
\end{rk}
If a tree $T$ of order $n$ satisfies $\alpha'(T)\geq k-1$, then by Theorem \ref{tree star}, $r(T,K_{1,k})=n$ provided $n \geq 2(3k-2)(2k-3)(k- 2) + 1$. Thus, combining Remark \ref{remark1}, we can see that Theorem \ref{tree star} implies Theorem \ref{tree}.

\vskip 2mm
A tree is a minimal connected graph. Theorem \ref{tree} tells us any graph $F$ of order $n$ with $\delta(F)\geq n-k$ contains all minimal connected spanning subgraphs with maximum degree not exceeding $n-2k+2$. 
A natural problem is, under the same condition on minimum degree, does $F$ contain some connected spanning subgraphs with more edges? To be precise, given a connected graph $G$ of order $n$ with at most $n(1+\epsilon)$ edges with bounded maximum degree, does $F$ contain $G$ as a spanning subgraph? Or when does $r(G,K_{1,k})=n$ hold? Here $\epsilon$ is a small positive constant depending on $k$ and such a graph is usually referred to as a sparse graph. In this paper, our main goal is to extend Theorem \ref{tree star} from trees to connected sparse graphs. 

The main result is the following.

\begin{thm}\label{mr}
   Let $k\geq 1$ and $n\geq 6k^3$. If $G$ is a connected graph with $n$ vertices, at most $n(1+1/(24k-12))$ edges and $\alpha'(G)=\alpha'$, then 
$$\max\{n,n+k-1-\alpha'-\beta \}\leq r(G,K_{1,k})\leq \max\{n,n+k-1- \alpha' \},$$
where $\beta=0$ if $k\,|\, n+k-2-\alpha'$, and $\beta=1$ otherwise. 
\end{thm}
For trees satisfying $n\geq 6k^3$, Theorem \ref{mr} gives the same Ramsey bounds as Theorem \ref{tree star}. Asymptotically, it lowers the sufficient order threshold from $12k^3+o(k^3)$ to $6k^3$.
Moreover, observe if $\alpha'(G)\geq k-1$, then $r(G,K_{1,k})=n$ in Theorem \ref{mr}, which is equivalent to saying that there exists a sparse connected spanning subgraph $G$ in a graph with given minimum degree, as stated below.

\begin{cor} \label{mr0}
 If $k\geq 1$ and $n\geq 6k^3$, then every graph $F$ of order $n$ with $\delta(F) \geq n -k$ contains every connected graph $G$ with $n$ vertices, at most $n(1+1/(24k-12))$ edges and $\alpha'(G)\geq k-1$.    
\end{cor}

By Remark \ref{remark1}, we have $\alpha'(T)\ge k-1$ for a tree $T$ if $\Delta(T)\le n-2k+2$. Consequently, for $n\geq 6k^3$, Corollary \ref{mr0} recovers the conclusion of Theorem \ref{tree} and gives a smaller asymptotic sufficient order threshold.
 
\begin{rk} \label{remark2} 
For a general graph $G$,  the parameter $\alpha'(G)$ is also closely related to $\Delta(G)$: if $\Delta(G)$ is large, then $\alpha'(G)$ would be small. 

The case $k=1$ is immediate, since the required inequality is $\alpha'(G)\geq0$. Hence assume $k\geq2$.
Let $G$ be a connected graph on $n$ vertices and at most $n(1+1/(24k-12))$ edges, where $n \geq 120k^2-180k+60$. If $\Delta(G)< n\left(1-1/(24k-12)\right)$, then for any vertex $v\in V(G)$,  $G_v$ has $m_v\geq \lfloor n/(24k-12) \rfloor>0$ vertices and 
$$e(G_v)\leq e(G)-(n-m_v)\leq m_v+\lfloor n/(24k-12) \rfloor.$$
This implies the average degree 
$$\bar{d}(G_v) \leq \frac{2e(G_v)}{m_v}\le 2\left(m_v+\left\lfloor \frac{n}{24k-12} \right \rfloor \right)\cdot \frac{1}{m_v}\leq 4.$$

\noindent 
By the Caro--Wei theorem \cite{Alon}, we get
$$\alpha(G_v)\geq \frac{m_v}{1+\bar{d}(G_v)}\geq \frac{m_v}{1+4}\geq \frac{1}{5} \cdot \left\lfloor \frac{n}{24k-12} \right \rfloor \geq k-1,$$
which means $\alpha'(G)\geq k-1$. 
\end{rk}

Combining Corollary \ref{mr0} and Remark \ref{remark2},
we have the following corollary, which tells us when there exists a sparse connected spanning subgraph with bounded maximum degree in a graph with given minimum degree.
\begin{cor}
If $k\geq 1$ and $n \geq \max\{6k^3,120k^2-180k+60\}$, then every graph $F$ of order $n$ with $\delta(F) \geq n -k$ contains every connected graph $G$ with $n$ vertices, at most $n(1+1/(24k-12))$ edges and  $\Delta(G) < n(1-1/(24k-12))$.
   
\end{cor}
 
Finally, as a generalization of Theorem \ref{mr}, we establish an upper bound for the Ramsey number $r(G,tK_{1,k})$ of a sparse graph $G$ versus $t$ copies of $K_{1,k}$. This bound is sharp when $\alpha'(G)\geq k-1$.

\begin{thm}\label{mr2}
    Let $k\geq 1$, $t\geq 1$, $g(k,t)=\max\{21tk-3k+6, 24k-12\}$ and $n\geq 28t^2k^3$. If $G$ is a connected graph with $n$ vertices, at most $n(1+1/g(k,t))$ edges and $\alpha'(G)=\alpha'$, then
$$r(G,tK_{1,k})\leq \max\{n,n+k-1- \alpha' \}+t-1.$$  
In particular, if $\alpha'\geq k-1$, then $r(G,tK_{1,k})=n+t-1$.
\end{thm}

When we consider embedding trees into a given graph, a very useful tool for doing this is ``Tree Trichotomy", developed by Burr, Erd\H{o}s, Faudree, Rousseau, and Schelp in \cite{Tree Trichotomy}, who showed that any
 tree must contain a long suspended path, i.e., a path in which all internal vertices have
 degree 2, or a large matching consisting of end-edges, or a large star consisting of end-edges, where an end-edge is one incident with an end-vertex, i.e., a vertex of degree 1. This method can help us embed a small tree first, and then extend the small tree to the tree as expected. However, when we consider embedding a sparse graph into a given graph, Tree Trichotomy does not work. Recently,
Zhang and Chen developed an enhanced version of the Tree Trichotomy Lemma for sparse graphs as below, which provided the clear structure of sparse graphs similar to that of trees. Therefore, this is the key ingredient for proving our results.
\begin{lem}(Zhang and Chen \cite{ZC}) \label{trichotomy}
Let $G$ be a connected graph with $n$ vertices and $n+\ell$ edges that is not a star, where $\ell \geq -1$, $n\geq q \geq 3$, and $s\geq2$ is an integer. If $G$ contains neither a suspended path of order $q$ nor a matching consisting of $s$ end-edges, then the number of vertices of degree at least 2 in $G$ is at most $\gamma$, and 
$G$ has a vertex adjacent to at least $\lceil \frac{n-\gamma}{s-1} \rceil$ vertices of degree 1, where $\gamma=(q-2)(2s+3\ell-2)+1$.
\end{lem}

 We conclude this section by introducing some additional notation and definitions. For a graph $G$, let $|G|$ and $e(G)$ denote the number of vertices and edges, respectively. For two graphs $G$ and $H$, $G\cup H$ denotes their disjoint union, and $G + H$ represents the disjoint union of two
 graphs $G$ and $H$, together with the new edges connecting every vertex of $G$  to every vertex of $H$. The
 graph $G-H$ will be the subgraph of $G$ induced by the vertices of $G$ not in $H$. Given a red-blue edge-coloring of $K_N$, $K_N[R]$ and $K_N[B]$ denote the spanning subgraphs of $K_N$ whose edge sets are the red edges and the blue edges, respectively.
 The red neighborhood of a vertex $v$, which is denoted by $N_R(v)$, is defined as the set of vertices that are adjacent to $v$ via
 red edges. For a complete graph $K_N$ and a vertex subset $A$, we use $K_N[A]$ to denote
 the subgraph induced by $A$. A cut vertex of a graph is a vertex whose removal (along with its incident edges) increases the number of connected components in the graph. A unicyclic graph is a connected graph containing exactly one cycle.

\section{Proof of Theorem \ref{mr}}
We divide the proof of Theorem \ref{mr} into two parts: the lower bound and the upper bound.
\subsection{The lower bound} 
On one hand, because the graph $K_{n-1}$ has no connected subgraph
with $n$ vertices and its complement has no $K_{1,k}$ for $k\geq 1$, it follows that $r(G,K_{1,k})\geq n$. 

On the other hand, it is easy to see that $n+k-1-\alpha'-\beta \ge n$ only if $\alpha'\le k-1$, so we need only to show  $r(G,K_{1,k})\geq n+k-1-\alpha'-\beta$ under the assumption that $\alpha' \le k-1$ and $k\geq 2$. 
Let $u$ be the vertex with $\alpha(G_u)=\alpha'$. Set
$n+k-2-\alpha'-\beta=tk+s$, where $0<s\leq k$. Since $n\geq 6k^3$, we have 
$$n+k-2-\alpha'-\beta=(t-k +1+ s)k + (k - s)(k - 1).$$
Let $F$ be the graph whose complement is the graph 
$$(t-k+1+s)K_k\cup (k-s)K_{k-1}.$$
Obviously, $\overline{F}$ contains no $K_{1,k}$. We now show that $F$ does not contain $G$ as a subgraph. Suppose to the contrary that $G$ is a subgraph of $F$. Recall the vertex $u\in V(G)$ satisfying $\alpha(G_u)=\alpha'$. Clearly, $u$ lies in some independent set $S$ of $F$ with at least $k-\beta$ vertices. Because all vertices of $G$ in $S$ must form an independent set of $G$ containing $u$, $S$ contains at most $\alpha' + 1$ vertices of $G$. Therefore, $F$ must have at least
$$n+k-\beta-(\alpha' + 1)=n+k-1-\alpha'-\beta$$
vertices, a contradiction. Thus we have $r(G,K_{1,k})\geq n+k-1-\alpha'-\beta$, and so the lower bound follows as expected.

\subsection{The upper bound} 
\noindent
The main method for establishing the upper bound is as follows. By Lemma \ref{trichotomy}, we can obtain a sparse $H$ of small order from $G$ through shortening a suspended path or deleting some vertices of degree 1. We first show that there is a red $H$, and then extend this $H$ to a red $G$. To do this, we need the following technical lemmas.

\begin{lem}(Hall \cite{Hall}) \label{hall1}
A bipartite graph $G=G[X,Y]$ has a matching which covers every vertex in $X$ if and only if $|N(S)|\geq |S|$ for all $S\subseteq X$.
\end{lem}

\begin{lem}(Chv\'atal~\cite{Chv})\label{Chv} For a tree $T$ on $n$ vertices,
	$r(T,K_m)=(n-1)(m-1)+1\,.$
\end{lem}
\begin{lem}(Burr~\cite{Burr})\label{Burr}
  For a tree $T$ on $n$ vertices, we have
  $r(T,K_{1,k})\le n+k-1$.
\end{lem}

Let $G-u$ denote the graph obtained from $G$ by removing the vertex $u$ and all edges incident with $u$.
\begin{lem}\label{2}
  For any graph $G$ and $u\in V(G)$, we have
  $r(G,K_{1,k})\leq r(G-u,K_{1,k})+k$.
\end{lem}
\noindent
{\bf Proof.}
Let $N=r(G-u,K_{1,k})+k$. Consider any red-blue edge-coloring of $K_N$ such that there is no blue $K_{1,k}$. Then for every $w\in V(K_N)$, we have
$$|N_R(w)|\geq N-k=r(G-u,K_{1,k}).$$
Therefore, the red neighborhood $N_R(w)$ contains a red copy of $G-u$. Since $w$ is red-adjacent to every vertex in this copy of $G-u$, we can obtain a red copy of $G$ by adding $w$ and identifying it with the vertex $u$.\qed
\vskip 2mm
For $uv\in E(G)$, let $G-uv$ denote the graph obtained from $G$ by removing the edge $uv$. It is clear that $G-u$ is a subgraph of $G-uv$, and thus we have the following immediate corollary.

\begin{cor}\label{cor}
  For a graph $G$ and $uv\in E(G)$, we have
  $r(G,K_{1,k})\leq r(G-uv,K_{1,k})+k$.
\end{cor}

\begin{lem}\label{lem:unicyclic}
  For a unicyclic graph $G$ on $n$ vertices, we have
  $r(G,K_{1,k})\le n+2k-2$.
\end{lem}

\noindent
{\bf Proof.}
Let $u$ be a vertex on the unique cycle of $G$. Then the graph $G-u$ is acyclic and has $n-1$ vertices, so it must be a subgraph of some tree $T$ on $n-1$ vertices. By Lemma~\ref{Burr},
$r(G-u,K_{1,k})\le r(T,K_{1,k})\le (n-1)+k-1=n+k-2$.
Then, by Lemma~\ref{2}, $r(G,K_{1,k})\le r(G-u,K_{1,k})+k\le n+2k-2$.
The proof is complete. \qed

\begin{lem}\label{lem:smalldegree}
  Let $n$ and $k$ be positive integers with $n\ge 2k+1$. If $G$ is a connected graph with $n$ vertices and at most $n(1+1/(2k+1))$ edges, then $G$ contains at least $k$ vertices with degree at most $2$.
\end{lem}

\noindent {\bf Proof.} Suppose to the contrary that $G$ contains at most $k-1$ vertices with degree at most $2$. Then, the number of vertices in $G$ with degree at least $3$ is at least $n-k+1$. Since $G$ has no isolated vertices, its total degree is at least $3(n-k+1)+(k-1)=3n-2k+2$.
On the other hand, since the number of edges in $G$ is at most $n(1+1/(2k+1))$, we have
$2n(1+1/(2k+1))\ge 2e(G)\ge 3n-2k+2$.
However, it follows from $n\ge 2k+1$ that $2n(1+1/(2k+1))<3n-2k+2$, which contradicts the previous inequality. \qed
\vskip 2mm

\noindent
The following lemma will be used several times to find a red subgraph $H$ of $G$.

\begin{lem}\label{3}
   Let $n$ and $k$ be positive integers. If $G$ is a connected graph with $n$ vertices and at most $n(1+1/(2k+1))$ edges, then
  \[r(G,K_{1,k})\le n+2k-2.\]
\end{lem}

\noindent {\bf Proof.} The result is trivial for $k=1$. Thus, we assume $k\ge 2$ and proceed by induction on the number of vertices in $G$.

When $n\le 2k$, the number of edges in the graph $G$ is at most $n(1+1/(2k+1))\le n+2k/(2k+1)$, which implies that $G$ has at most $n$ edges. Since $G$ is connected, it must be either a tree or a unicyclic graph. By Lemmas~\ref{Burr} and~\ref{lem:unicyclic}, the lemma follows. Next, we assume $n\ge 2k+1$.

According to Lemma~\ref{lem:smalldegree}, $G$ contains at least $k$ vertices with degree at most $2$. Now we perform a ``reduction'' operation on certain vertices in $G$ with degree at most $2$. Specifically, for a vertex of degree $1$ or a non-cut vertex of degree $2$, the reduction operation refers to deleting this vertex and its incident edges. For a cut vertex of degree $2$, the reduction operation refers to deleting this vertex and its incident edges, and then adding an edge between its two neighbors.

    In the first step, we arbitrarily select a vertex with degree at most $2$ in $G$ and perform the reduction operation, obtaining a new graph denoted as $G_1$. For $i\in [k-1]$, in step $i+1$, we arbitrarily select a vertex with degree at most $2$ in $G_i$ and perform the reduction operation, obtaining a new graph denoted as $G_{i+1}$. Each reduction preserves connectedness and does not increase the degree of any surviving vertex: deleting a vertex of degree $1$ or a non-cut vertex of degree $2$ can only decrease degrees, while suppressing a cut vertex of degree $2$ preserves the degrees of its two neighbors. Hence the original set of $k$ vertices of degree at most $2$ always leaves an eligible vertex for the next step. The operation can therefore be performed $k$ times. Each step deletes one vertex and decreases the number of edges by at least one, so $G_k$ is connected, $|G_k|=|G|-k$, and $e(G_k)\le e(G)-k$.
    
    Let $N=n+2k-2$. We aim to prove that for any red-blue edge-coloring of the complete graph $K_N$, either there exists a red subgraph isomorphic to $G$, or there exists a blue subgraph isomorphic to $K_{1,k}$. Assuming the latter does not exist, it suffices to show that $K_N$ contains a red subgraph isomorphic to $G$.

    We first prove that $K_N$ contains $G_k$ as a red subgraph. It is easy to see that
    \[e(G_k)\le n\left(1+\frac{1}{2k+1}\right)-k\le (n-k)\left(1+\frac{1}{2k+1}\right)+1.\]
    If
    $$e(G_k)\le (n-k)\left(1+\frac{1}{2k+1}\right),$$
    then by the induction hypothesis, $r(G_k,K_{1,k})\le (n-k)+2k-2<N$.
    Therefore, $K_N$ contains a red subgraph isomorphic to $G_k$. If
    $$e(G_k)=\left\lfloor (n-k)\left(1+\frac{1}{2k+1}\right)\right\rfloor+1,$$
    then there exists an edge $e$ in $G_k$ such that $G_k-e$ remains connected. By Corollary~\ref{cor} and the induction hypothesis,
    \[r(G_k,K_{1,k})\le r(G_k-e,K_{1,k})+k\le (n-k)+2k-2+k=N.\]
    Thus, $K_N$ still contains a red subgraph isomorphic to $G_k$.

    Recall that the graph $G_k$ is obtained from $G$ through $k$ successive reduction operations. Next, we attempt to reconstruct $G$ in $K_N[R]$ in a reverse manner. That is, for each $i\in [k]$, if $K_N$ contains $G_i$ as a red subgraph, then it also contains $G_{i-1}$ as a red subgraph, where $G_0=G$.

    Suppose that at some step the reconstruction fails, meaning that there exists some $i\in [k]$ such that $K_N$ contains $G_i$ as a red subgraph but does not contain $G_{i-1}$ as a red subgraph. Without loss of generality, assume that $G_i$ is obtained from $G_{i-1}$ by reducing a vertex $u_0$. Additionally, let $U_i$ denote the set of vertices in $K_N$ that are not in $V(G_i)$. Then,
    \[|U_i|=N-|G_i|\ge n+2k-2-(n-1)= 2k-1.\]

    If $u_0$ is a vertex of degree $1$ in $G_{i-1}$, let its neighbor in $G_{i-1}$ be denoted as $u_1$. Since $u_1$ is present in $G_i$, we may assume that $u_1$ is embedded in $K_N$ as vertex $v_1$. Because $K_N$ does not contain a blue $K_{1,k}$, the number of blue edges between $v_1$ and $U_i$ is at most $k-1$. Consequently, there must be at least one red edge between $v_1$ and $U_i$. Suppose $v_1$ is red-adjacent to some $w\in U_i$. By embedding $u_0$ at vertex $w$, we can find a red subgraph isomorphic to $G_{i-1}$ in $K_N$, leading to a contradiction.

    If $u_0$ is a vertex of degree $2$ in $G_{i-1}$, let its neighbors in $G_{i-1}$ be $u_1$ and $u_2$. Since both $u_1$ and $u_2$ are present in $G_i$, we may assume that they are embedded in $K_N$ as vertices $v_1$ and $v_2$, respectively. Because $K_N$ does not contain a blue $K_{1,k}$, both $v_1$ and $v_2$ have at most $k-1$ blue edges connecting them to $U_i$.
    Thus there is at least $2k-1-2(k-1)= 1$ vertex in $U_i$ must be red-adjacent to both $v_1$ and $v_2$. Let this vertex be $w$. By embedding $u_0$ at vertex $w$, we can find a red subgraph isomorphic to $G_{i-1}$ in $K_N$, again a contradiction. This completes the proof of the lemma.   \qed
\vskip 2mm    

\vskip 5mm
Now, we are in a position to establish the upper bound.

Let $N=\max\{n,n+k-1- \alpha' \}$. Given a red-blue edge-coloring of $K_N$, assume that there is no blue $K_{1,k}$. We will show that $K_N$ contains a red $G$. The result is trivial for $k = 1$, so assume $k \geq2$.

We consider the following three cases separately.
\vskip 2mm
\noindent {\textbf{Case 1.}} $G$ has a suspended path with at least $4k-2$ vertices.
\vskip 2mm

Let $H$ be the graph on $n-(2k-2)$ vertices obtained from $G$ by shortening the suspended path by $2k-2$ vertices. Clearly, by $n\geq 6k^3$,
\begin{align*}
e(H)& \leq  n \left(1+\frac{1}{24k-12} \right)-(2k-2) \leq  (n-(2k-2))\left(1+\frac{1}{2k+1}\right).
\end{align*}
Then applying Lemma \ref{3}, $r(H,K_{1,k})\leq n-(2k-2)+2k-2=n\leq N$. Thus $K_N$ contains a red $H$. Of course $G$ can be obtained
from $H$ by lengthening the suspended path in $H$ by $2k-2$ vertices. Let $H'$ be a
subgraph of $K_N[R]$ in which this suspended path has been lengthened as much as
possible (up to $2k-2$). If $H'$ is isomorphic to $G$, the proof of this case is complete.
If not, let $P$ be the 
suspended path in $H'$. By the maximality of $P$, 
$v$ is not red-adjacent to two consecutive vertices of $P$ for $v\in V(K_N)-V(H')$. Since $|P|\geq (4k-2)-(2k-2)=2k$, $v$ has at least $k$ blue neighbors in $P$, which implies that $K_N$ contains a blue $K_{1,k}$, a contradiction.

\vskip 2mm
 \noindent {\textbf{Case 2.}} 
$G$ has a matching consisting of at least $2k-2$ end-edges.
\vskip 2mm

Let $H$ be the graph on $n-(2k-2)$ vertices obtained from $G$ by deleting $2k-2$ end-vertices of this matching. Because $n\geq 6k^3$, 
\begin{align*}
e(H)& \leq  n\left(1+\frac{1}{24k-12}\right)-(2k-2) \leq  (n-(2k-2))\left(1+\frac{1}{2k+1}\right).
\end{align*}
Then by Lemma \ref{3}, 
$$r(H,K_{1,k})\leq n-(2k-2)+2k-2=n\leq N.$$ 
Thus $K_N$ contains a red $H$. 
Let $X$ be the set of $2k-2$ vertices of $H$ adjacent to the end-vertices of $G$ deleted.  Select a set $Y$ from $V(K_N)-V(H)$ with $|Y|=2k-2$. If there is a red matching in $K_N$ between $X$ and $Y$ which saturates $X$, then $K_N$
contains a red $G$. If not, then by Lemma \ref{hall1}, there is a nonempty subset $X'\subseteq X$, such that $Y' =N_R[X'] \cap Y$ satisfies $|Y'|<|X'|$. Because $K_N$ contains no blue $K_{1,k}$, each vertex of $X$ is blue-adjacent to at most $k-1$ vertices of $Y$. Thus $|Y'|\geq 2k-2-k+1=k-1$ and then $|X'|\geq k$. Note that any vertex of $Y-Y'$ is blue-adjacent to all vertices of $X'$, which gives a blue $K_{1,k}$ in $K_N$, a contradiction.

\vskip 2mm
\noindent {\textbf{Case 3.}} 
$G$ has neither a suspended path of $4k-2$ vertices, nor a matching formed by $2k-2$ end-edges.
\vskip 2mm
To apply Lemma \ref{trichotomy} in a more straightforward manner, we let $e(G)=n+\ell$, $q=4k-2$ and $s=2k-2$, where $\ell\leq n/(24k-12)$.

If $G$ is not a star, then by Lemma \ref{trichotomy},  $G$ has a vertex adjacent to at least 
$$\left\lceil \frac{n-\gamma}{2k-3} \right\rceil$$ 
vertices of degree 1, where $$\gamma=(q-2)(2s+3\ell-2)+1 \leq (4k-4)\left(4k-6+ \frac{3n}{24k-12}\right)+1.$$ 

If $G$ is a star, then $G$ has a vertex adjacent to $n-1$ vertices of degree 1.
Denote this vertex by $v$. Since
$\gamma\le 16k^2-40k+25+\frac{k-1}{2k-1}n$,
the assumptions $n\geq6k^3$ and $k\geq2$ give
\[
n-\gamma-(2k-3)(k^2-1)
\geq
\frac{2k^4-24k^3+97k^2-98k+28}{2k-1}>0.
\]
Consequently,
\[
\min \left\{ \left\lceil \frac{n-\gamma}{2k-3} \right\rceil, \, n-1 \right\} \ge k^2.
\]

Let $H$ be the graph on $n-k^2$ vertices obtained from $G$ by deleting $k^2$ vertices of degree 1, which 
are adjacent to $v$. Because $n\geq 6k^3$,
\begin{align*}
e(H)& \leq  n\left(1+\frac{1}{24k-12}\right)- k^2 \leq  (n- k^2)\left(1+\frac{1}{2k+1}\right). \end{align*}
By Lemma \ref{3}, we have 
\begin{align*}
r(H,K_{1,k}) \leq n-k^2+2k-2 \leq N-(k^2-2k+2).
\end{align*}
Given $u\in V(K_N)$, select a maximal set $A\subseteq N_R(u)$ such that each vertex of $N_B(u)$ is red-adjacent to each vertex of $A$.
Since $K_N$ contains no blue $K_{1,k}$, each vertex of $N_B(u)$ is blue-adjacent to at most $k-2$ vertices of $N_R(u)$ and $|N_B(u)| \leq k-1$. Thus 
\begin{align*}
|A| & \geq |N_R(u)|-|N_B(u)|(k-2)\\
& \geq N-k-(k-1)(k-2)\\
& =N-(k^2-2k+2).
\end{align*}
By the arguments above, we can see that $K_N[A]$ contains a red $H$. If $|N_B(u)|=0$, then $|N_R(u)|=N-1\geq n-1$. After replacing the vertex $v$ of this $H$ by $u$, at least $(n-1)-(|H|-1)=k^2$ unused red neighbors of $u$ remain, and they accommodate the deleted end-vertices. This gives a red copy of $G$ in $K_N$. If $|N_B(u)|>0$, then replacing the vertex $v$ of this $H$ by $u$, we get a new red copy $H'$ of $H$.
Set $H'_u=H'-N_{H'}[u]$ and $H_v=H-N_H[v]$. Select an independent set $I$ of $\min \{\alpha(H_u'), |N_B(u)|\}$ vertices in $H_u'$. Replace the independent set $I$ of $H'_u$
 with $|I|$ vertices of $N_B(u)$. This gives another red copy $H''$ of $H$ because each vertex of $N_B(u)$ is red-adjacent to each vertex of $A$. We claim that 
 $$ |N_R(u)-V(H'')| \geq |N_R(u)|+|I|+1-|H|\geq  k^2.$$
 In fact, if $\alpha(H_u')\geq |N_B(u)|$, then $|I|=|N_B(u)| \leq k-1$. Thus 
\begin{align*}
& |N_R(u)-V(H'')|\\
\geq~& |N_R(u)|+ |N_B(u)|+1-|H|\\
=~& N-|H|\\
\geq~& n-(n- k^2)\\
=~& k^2.
\end{align*}
If $\alpha(H_u')\leq |N_B(u)|$, then $|I|=\alpha(H_u')$. 
 Thus
 \begin{align*}
& |N_R(u)-V(H'')| \\
  \geq~& |N_R(u)|+\alpha(H_u')+1-|H|\\
 \geq~ & N-k+\alpha(H_u')+1-n+ k^2\\
 \geq~ & n+k-1-\alpha'-k+\alpha(H_u')+1-n+ k^2\\
 \geq~ &  k^2.
 \end{align*} 
The last inequality holds because $\alpha(H_u')=\alpha(H_v)=\alpha(G_v)\geq \alpha'$.
Therefore, the vertex $u$ is red-adjacent to at least $k^2$ vertices in $V(K_N)\backslash V(H'')$, which implies that $K_N$ contains a red $G$. This completes the proof of the theorem.\qed

\section{Proof of Theorem \ref{mr2}}
\noindent
The main idea for proving Theorem \ref{mr2} is by induction on $t$, together with the main method used in Theorem \ref{mr}. Moreover, we need the following two additional lemmas.
\begin{lem} \label{tSkweak}
   Let $G$ be a connected graph with $n$ vertices, at most $n(1+1/(24k-12))$  edges and $\alpha'(G)=\alpha'$. Then for $k\geq 1$, $t\geq 1$ and $n\geq 6k^3$,
$$r(G,tK_{1,k})\leq \max\{n,n+k-1- \alpha' \}+(t-1)(k+1).$$  
\end{lem}

\noindent
{\bf Proof.} Consider an arbitrary red-blue edge-coloring of a complete graph on 
$\max\{n,n+k-1- \alpha' \}+(t-1)(k+1)$
vertices. We begin by finding as many disjoint blue copies of $K_{1,k}$ as possible within the graph. Suppose we find $s$ such blue copies of $K_{1,k}$. If $s \geq t$, then the proof is complete. 
Thus, assume that $s \leq t-1$. There are at least $ \max \{n,n+k-1- \alpha' \}$ vertices outside the blue $s K_{1,k}$.
Since the complete graph induced by these vertices does not contain a blue $K_{1,k}$, by Theorem \ref{mr}, there must exist a red $G$. \qed

\begin{lem} (Hall \cite{Hall}) \label{matching}
 Consider a complete bipartite graph $K_{a,b}$, where $a \leq b$, with parts $X = \{x_1,\ldots,x_a\}$ and $Y=\{y_1,y_2,\ldots,y_b\}$, whose edges are colored red and blue. Then, one of the following holds:
 
 \noindent
 (1) There exists a red matching of size $a$;

 \noindent
 (2) For some $0 \leq c \leq a-1$, there exists a blue subgraph $K_{c+1,b-c}$, where $c+1$ vertices are in $X$.   
 \end{lem}
\vskip 3mm
We now begin to  prove Theorem \ref{mr2}. 

We first dispose of the case $k=1$. Then $tK_{1,1}=tK_2$. For $t=1$, the assertion is immediate. For $t=2$, Zhang and Chen \cite{ZC} proved that $r(G,2K_2)=n+1$ for every connected noncomplete graph $G$ of order $n$; our graph is noncomplete under the present size condition. For $t\geq3$, they proved that $r(G,tK_2)=n+t-1$
whenever $n\geq(4t-5)(2t-3)+3$ and $e(G)\leq n+n^2/(4t-5)-2$. Both inequalities follow from $n\geq28t^2$ and $e(G)\leq n(1+1/g(1,t))$. Hence the theorem holds when $k=1$.

We may therefore assume $k\geq2$. The proof is by induction on $t$. When $t=1$, we have $g(k,1)\geq24k-12$, so the assertion follows from Theorem \ref{mr}. Assume that the theorem holds for $t-1$, and let $t\geq2$. Notice that $g(k,t)=21tk-3k+6$ for $t\geq2$.

 To apply the inductive method, it is necessary to ensure that the lower bound on $n$ and the upper bound on $e(G)$ are
 compatible with induction. This can be derived from the facts that the lower bound on $n$ is increasing in $t$ and the upper bound on $e(G)$ is decreasing in $t$.
 
 Let $N=\max\{n,n+k-1- \alpha' \}+t-1$. Given a red-blue edge-coloring of $K_N$, assume that $K_N$ contains neither a red $G$ nor a blue $tK_{1,k}$.

\vskip 2mm
\noindent {\textbf{Case 1.}} $G$ has a suspended path with at least $(3t-1)k$ vertices.
\vskip 2mm

Let $H$ be the graph on $n-(t-1)k$ vertices obtained from $G$ by shortening the suspended path by $(t-1)k$ vertices. Because $n\geq 28t^2k^3$ and $t\geq 2$,
\begin{align*}
|H| &=n-(t-1)k \geq 6k^3,\\
e(H)& \leq  n\left(1+\frac{1}{21tk-3k+6}\right)-(t-1)k \leq  (n-(t-1)k)\left(1+\frac{1}{24k-12}\right).
\end{align*}
The remaining suspended path in $H$ has at least $2tk$ vertices. For every $z\in V(H)$, deleting $N_H[z]$ removes at most three vertices from this path. Hence
\[
\alpha(H-N_H[z])\geq \left\lceil\frac{2tk-3}{2}\right\rceil=tk-1\geq k-1.
\]
Thus $\alpha'(H)\geq k-1$. By Lemma \ref{tSkweak},
\[
r(H,tK_{1,k})\leq |H|+(t-1)(k+1)=n+t-1\leq N,
\]
which implies that $K_N$ contains a red $H$. Let $H'$ be a
subgraph of $K_N[R]$ in which this suspended path has been lengthened as much as
possible (up to $(t-1)k$). If $H'$ is isomorphic to $G$, the proof of this case is complete.
If not, let $P$ be the 
suspended path in $H'$. By the maximality of $P$, $v$ is not red-adjacent to two consecutive vertices of $P$ for any $v\in V(K_N)-V(H')$. Since $|P|\geq (3t-1)k-(t-1)k=2tk$, each vertex of $V(K_N)-V(H')$ has at least $tk$ blue neighbors in $P$. Moreover, $|V(K_N)-V(H')| \geq N-(n-1)\geq n+t-1-(n-1)=t$. Thus we can greedily find a blue $tK_{1,k}$ in $K_N$, a contradiction.

\vskip 2mm
\noindent {\textbf{Case 2.}} 
$G$ has a matching consisting of at least $2tk-2$ end-edges.
\vskip 2mm
By the induction hypothesis, we have
 $$r(G,(t-1)K_{1,k})\leq \max\{n,n+k-1- \alpha' \}+t-2<N,$$
 which implies that the graph $K_N$ contains a blue subgraph $(t-1)K_{1,k}$. We denote the
 vertex set of this subgraph by $A$.

Let $H$ be the graph on $n-2tk+2$ vertices obtained from $G$ by deleting $2tk-2$ end-vertices of this matching. Because $n\geq 28t^2k^3$ and $t\geq 2$,
\begin{align*}
|H| &=n-2tk+2 \geq  6k^3,\\
e(H) & \leq  n\left(1+\frac{1}{21tk-3k+6}\right)-2tk+2 \leq  (n-2tk+2)\left(1+\frac{1}{24k-12}\right).
\end{align*}
By Lemma \ref{3},
\[
r(H,K_{1,k})\leq |H|+2k-2=n-2k(t-1)\leq N-|A|,
\]
which implies that $K_N-A$ contains a red $H$.
Let $X$ be the set of $2tk-2$ vertices of $H$ adjacent to the end-vertices of $G$ deleted and $Y=V(K_N)-V(H)$. By Lemma \ref{matching}, there is either a red matching that covers $X$, or a blue $K_{c+1,|Y|-c}$, where $0 \leq c \leq |X|-1$. If the former holds, then $K_N$ would contain a red $G$, leading to a contradiction. If the latter holds, since $K_N-A$ contains no blue $K_{1,k}$, it follows that $|Y|-c\leq |A|+k-1$. Then
\begin{align*}
   c&\geq |Y|-|A|-k+1 \\
   &=N-|H|-(t-1)(k+1)-k+1 \\
   &= N-n+2tk-2-(t-1)(k+1)-k+1 \\
   &\geq t-1+2tk-2-tk-t+k+1-k+1 \\
   &= tk-1.
\end{align*}
Note that $|Y|-c\geq |Y|-(|X|-1)=N-(n-|X|)-(|X|-1)\geq t$.
 Thus $K_{c+1,|Y|-c}$ contains $K_{tk,t}$ as a subgraph, which implies that $K_N$ contains a blue $tK_{1,k}$, a contradiction.

\vskip 2mm
\noindent {\textbf{Case 3.}} 
$G$ has neither a suspended path of $(3t-1)k$ vertices, nor a matching formed by $2tk-2$ end-edges.
\vskip 2mm
To apply Lemma \ref{trichotomy} in a more straightforward manner, we let $e(G)=n+\ell$, $q=(3t-1)k$ and $s=2tk-2$, where $\ell\leq n/(21tk-3k+6)$.

If $G$ is not a star, then by Lemma \ref{trichotomy},  $G$ has a vertex $v$ adjacent to at least 
$$\left\lceil \frac{n-\gamma}{2tk-3} \right\rceil$$ 
vertices of degree 1, where 
$$\gamma=(q-2)(2s+3\ell-2)+1 \leq ((3t-1)k-2)\left(4tk-6+\frac{3n}{21tk-3k+6}\right)+1.$$ 
If $G$ is a star, then $G$ has a vertex $v$ adjacent to $n-1$ vertices of degree 1.
Since $21tk-3k+6\geq18tk$, the displayed upper bound gives
$\gamma\leq 12t^2k^2+ n/2+1$.
Using $n\geq28t^2k^3$ and $k,t\geq2$, we obtain
\begin{align*}
n-\gamma
\geq \frac n2-12t^2k^2-1\geq 14t^2k^3-12t^2k^2-1>(2tk-3)\bigl((t-1)k+k^2-k+1\bigr).
\end{align*}
The same conclusion is immediate for the star, because its center has $n-1$ end-neighbors. Therefore,
$$\min \left \{\left\lceil \frac{n-\gamma}{2tk-3} \right\rceil,n-1 \right \}  \geq (t-1)k+k^2-k+2.$$
By the induction hypothesis, we have
 $$r(G,(t-1)K_{1,k})\leq \max\{n,n+k-1- \alpha' \}+t-2<N,$$
 which implies that the graph $K_N$ contains a blue subgraph $(t- 1)K_{1,k}$. We denote the
 vertex set of this subgraph by $A$.

 First, we claim that there is a vertex $x$ in $K_N-A$, which is red-adjacent to at least $(t-1)k$ vertices of $A$. Assume not, then each vertex in $K_N-A$ is blue-adjacent to at least $t$ vertices of $A$. Thus the number of blue edges between $A$ and $V(K_N)-A$ is at least $(N-|A|)t$.
 If there are two vertices in some $K_{1,k}$ contained in $A$ which each receive at least $2k$ blue edges from $V(K_N)-A$, then we may choose two disjoint sets of $k$ blue neighbors and replace this $K_{1,k}$ by $2K_{1,k}$. Together with the remaining blue $(t-2)K_{1,k}$, this gives a blue $tK_{1,k}$. If not, then the number of blue edges between $V(K_N)-A$ and $A$ is at most $(N-|A|+(2k-1)k)(t-1)$. Therefore,
 $$(N-|A|)t\leq (N-|A|+(2k-1)k)(t-1),$$
 which is equivalent to
 $N-|A|\leq (2k-1)k(t-1)$,
 that is,
 $N\leq (2k^2+1)(t-1)$.
 This is impossible because $N>n\geq 28t^2k^3$, and hence the claim is true.

 Let $d=(t-1)k+k^2-2k+2$, and let $H$ be the graph obtained from $G$ by deleting exactly $d$ end-vertices adjacent to $v$. The preceding estimate ensures that at least $k$ end-vertices adjacent to $v$ remain in $H$. Consequently,
 \[
 \alpha'(H)\geq\min\{\alpha'(G),k-1\}.
 \]
 Indeed, for $z=v$ we have $H-N_H[v]=G-N_G[v]$, while for $z\neq v$ at least $k-1$ of the retained end-vertices belong to $H-N_H[z]$. Because $n\geq 28t^2k^3$ and $t\geq 2$,
\begin{align*}
|H|& =n-((t-1)k+k^2-2k+2) \geq  6k^3,\\
    e(H) &\leq n\left(1+\frac{1}{21tk-3k+6}\right)-((t-1)k+k^2-2k+2)\\ 
    & \leq  (n-((t-1)k+k^2-2k+2))\left(1+\frac{1}{24k-12}\right).
\end{align*}
Applying Theorem \ref{mr}, we have
\begin{align*}
r(H,K_{1,k})
&\leq \max\{|H|,|H|+k-1-\alpha'(H)\}\\
&\leq \max\{n,n+k-1-\alpha'\}-d\\
&=N-|A|-(k^2-2k+2).
\end{align*}
  For convenience, we set $V(K_N)-A=S$. Select a maximal set $U\subseteq N_R(x)\cap S$ such that each vertex of $N_B(x)\cap S$ is red-adjacent to each vertex of $U$.
Since $K_N-A$ contains no blue $K_{1,k}$, each vertex of $N_B(x)\cap S$ is blue-adjacent to at most $k-2$ vertices of $N_R(x)\cap S$ and $|N_B(x)\cap S| \leq k-1$. Thus 
\begin{align*}
|U| & \geq |N_R(x)\cap S|-|N_B(x)\cap S|(k-2)\\
& \geq |S|-k-(k-1)(k-2)\\
& =N-|A|-(k^2-2k+2).
\end{align*}
By the arguments above, we can see that $K_N[U]$ contains a red $H$.

If $|N_B(x)\cap S|=0$, then
 $$|N_R(x)\cap S|+|N_R(x)\cap A|\geq |S|-1+(t-1)k=N-(t-1)(k+1)-1+(t-1)k\geq n-1$$
and, after replacing the vertex $v$ of this $H$ by $x$, at least
$(n-1)-(|H|-1)=n-|H|=d$
unused red neighbors of $x$ remain. They accommodate exactly the $d$ deleted end-vertices and give a red copy of $G$.

If $|N_B(x)\cap S|>0$, then replacing the vertex $v$ of this $H$ by $x$, we get a new red copy $H'$ of $H$. Set $H'_x=H'-N_{H'}[x]$ and $H_v=H-N_H[v]$. Select an independent set $I$ of $\min \{\alpha(H_x'), |N_B(x)\cap S|\}$ vertices in $H_x'$. Replace the independent set $I$ of $H'_x$
 with $|I|$ vertices of $N_B(x)\cap S$. This gives another red copy $H''$ of $H$ because each vertex of $N_B(x)\cap S$ is red-adjacent to each vertex of $U$. We now verify the required number of unused red neighbors in $S$. Put $b=|N_B(x)\cap S|$. If $\alpha(H_x')\geq b$, then $|I|=b$ and
\begin{align*}
|(N_R(x)\cap S)\setminus V(H'')|
&=|N_R(x)\cap S|+b+1-|H|\\
&=|S|-|H|\\
&=\max\{n,n+k-1-\alpha'\}-n+k^2-2k+2\\
&\geq k^2-2k+2.
\end{align*}
If $\alpha(H_x')\leq b$, then $|I|=\alpha(H_x')=\alpha(H_v)=\alpha(G_v)\geq\alpha'$. Since $|N_R(x)\cap S|\geq |S|-k$, we obtain
\begin{align*}
|(N_R(x)\cap S)\setminus V(H'')|
&\geq |S|-k+\alpha'+1-|H|\\
&=\max\{n,n+k-1-\alpha'\}-n+\alpha'-k+1+k^2-2k+2\\
&\geq k^2-2k+2.
\end{align*}
 Recall that $x$ is red-adjacent to at least $(t-1)k$ vertices of $A$, we can see that $K_N$ contains a red $G$. 
 
 Up to now, we have shown that $r(G,tK_{1,k})\leq \max\{n,n+k-1- \alpha' \}+t-1$. Moreover, if $G$ satisfies $\alpha'(G)\ge k-1$, then we
 have $r(G,tK_{1,k})\leq n+t-1$. The graph $K_{n-1}\cup K_{t-1}$ contains no connected graph of order $n$. Its complement is $K_{n-1,t-1}$, and every copy of $K_{1,k}$ in this bipartite graph uses at least one vertex from the part of order $t-1$; hence it cannot contain $t$ vertex-disjoint copies. Therefore $r(G,tK_{1,k})\ge n+t-1$, and so $r(G,tK_{1,k})= n+t-1$.
 
 This completes the proof of the theorem.
 \qed

\section{Concluding Remark}
The numerical hypotheses in Theorems \ref{mr} and \ref{mr2} are convenient sufficient conditions arising from the third cases of their proofs; they are not claimed to be best possible. For $k\geq2$, in the special case $e(G)\leq n$, the order condition in Theorem \ref{mr} can be improved directly. Indeed, in Case 3 we have $\ell\leq0$ and hence
$\gamma\leq(4k-4)(4k-6)+1=16k^2-40k+25$.
If $n>2k^3+13k^2-40k+25$,
then $(n-\gamma)/(2k-3)>k^2-1$, so the vertex supplied by Lemma \ref{trichotomy} is adjacent to at least $k^2$ end-vertices. Cases 1 and 2 remain valid because the corresponding reduced graph $H$ satisfies $e(H)\leq |H|$. Thus the conclusion of Theorem \ref{mr} holds for every connected graph with at most $n$ edges under this smaller order condition. In particular, it holds for every tree or unicyclic graph of order $n$ satisfying this condition.

Substantial further improvements in the general edge and order bounds appear to require different methods.

 \section*{Acknowledgments}
This research was supported  by National Key R\&D Program of China under grant number 2024YFA1013900 and NSFC under grant number 12471327.

\end{document}